\documentclass[12pt,a4paper]{article}
\usepackage[T2A]{fontenc}
\usepackage[cp1251]{inputenc}
\usepackage[russian]{babel}
\usepackage{amsmath,amssymb,amsthm,amsfonts,amscd}
\begin{document}
\newtheorem{example}{Example}
\newtheorem{remark}{Remark}
\newtheorem{theorem}{Theorem}
\newtheorem{conjecture}{Conjecture}
\newtheorem{lemma}{Lemma}
\def\Z{{\Bbb Z}}
\def\H{{\Bbb H}}
\def\R{{\Bbb R}}
\def\RP{{\Bbb R}\!{\rm P}}
\def\N{{\Bbb N}}
\def\C{{\Bbb C}}
\def\Ca{{\Bbb Ca}}
\def\A{{\bf A}}
\def\D{{\bf D}}
\def\k{{\bf k}}
\def\E{{\bf E}}
\def\F{{\bf F}}
\def\V{\vec{\bf V}}
\def\L{{\bf L}}
\def\M{{\bf M}}
\def\c{{\bf c}}
\def\q{{\bf q}}
\def\i{{\bf i}}
\def\e{{\bf e}}
\def\f{{\bf f}}
\def\g{{\bf g}}
\def\h{{\bf h}}

\def\fr{{\operatorname{fr}}}
\def\st{{\operatorname{st}}}
\def\mod{{\operatorname{mod}\,}}
\def\cyl{{\operatorname{cyl}}}
\def\dist{{\operatorname{dist}}}
\def\grad{{\bf{grad}}}
\def\sign{{\operatorname{sign}}}
\def\rot{{\operatorname{rot}}}
\def\j{{\bf{j}}}
\def\R{{\Bbb R}}
\def\Z{{\bf{Z}}}
\def\z{{\bf{z}}}
\def\B{{\bf B}}
\def\e{{\bf e}}
\def\L{{\bf L}}
\def\valpha{\vec{\alpha}}
\def\vxi{\vec{\xi}}
\sloppy
\title{A remark on the space of 7-gons with a fixed total length in $\R^3$}
\author{Akhmet'ev P.M.}
\maketitle

\begin{abstract}
Based on the model of the space $Pol_3(n)$ of polygons in $R^3$ with limited number of vertex, which was proposed by Jean-Claude Hausmann and Allen Knutson, and developed by several authors: Jason Cantarella, Alexander Y. Grosberg, Robert Kusner, and Clayton Shonkwiler, we prove that there exists an isometric isotopy of $Pol_3(n)$, $n=7$, into itself, which
transforms an arbitrary polygon to its mirror copy, and, additionally, preserves lengths of projections of polygons into the two coordinate planes, and keeps projection of polygons onto the line. The proof is based on elementary arguments with Cayley numbers.
A possible generalization of the statement for greater $n$ is related with a theorem by I.James on strong Kervaire invariants in stable homotopy of spheres.
\end{abstract}
\[  \]

\section*{Introduction}

In the paper \cite{C-G-K-Sh} the symmetric measure on closed polygons of fixed total length is investigated.
The space $Pol_3(n)$  consists of vectors $(\vec{e}_1, \dots, \vec{e}_n) \in (\R^3)^n$, which satisfy the conditions:

--1.  $\sum_{i=1}^n \vec{e}_i = 0$.

--2. $\sum_{i=1}^n \vert \vec{e}_i \vert = 2$.

Such a vector $(\vec{e}_1, \dots, \vec{e}_n)$ represents a closed $n$-gone in $\R^3$ with a fixed (twice unite) total length. 
The authors apply an approach by Jean-Claude Hausmann and Allen Knutson to describe a parametrization of the space  $Pol_3(n)$
by the homogeneous Stiefel manifold $V_2(\C,n)$. This allow to introduce the symmetric measure on $Pol_3(n)$ using the quaternions.

In the present paper we use elementary properties of the Cayley numbers to prove the following theorem.

\begin{theorem}\label{main}
There exists a one-parameter isometric isotopy $\Theta$ of $Pol_3(n)$, $n=7$, to itself, which joins each closed polygon of the length $2$ with its mirror image. Additionally: 

--I this  isotopy preserves lengths of the projections of polygons into the two coordinate planes  $(x,y,0) \subset \R^3$, $(x,0,z) \subset \R^3$, 

--II this  isotopy keeps the projection of polygons onto the line $(x,0,0) \subset \R^3$.
\end{theorem}

A generalization of Theorem \ref{main} for spaces $Pol_3(n)$ with larger $n \ge 8$ could be interesting. I conjecture, such a generalization is possible for
$n=15$, using an application  by I.James of a Toda's result to the strong Kervaire Invariant One Conjecture, see \cite{J}. The theorem by James is proved using homotopy theory and explicit formulas are unknown. I am not able to formulate a conjecture precisely, because I am not able to formulate analogs of properties I, II in the case $n=15$. 

A partial solution of the Kervaire Invariant One Problem in stable homotopy of spheres by M.A. Hill, M.J. Hopkins, and D.C. Ravenel (2010), and a particular solution of a weaker result by the author of the Kervaire Invariant One Problem (my approach in a final form is not published), which is based on alternative arguments, allows to formulate the following conjecture.

\begin{conjecture}\label{arf}
For an arbitrary $n$, $n \ge n_0$, where $n_0$ is sufficiently large,  Theorem $\ref{main}$ (even without additional properties --I, --II) is not valid.
\end{conjecture}


The paper is organized as following. All required constructions from the paper \cite{C-G-K-Sh} are collected if Section 1. In Section 2 new elementary computations are presented and Theorem is proved. This results was initiated by several discussions with the first and the forth authors of \cite{C-G-K-Sh}. The  results was prepared at Prof. D.D.Sokoloff's seminaire in Moscow, May (2013). 

\section{The model spaces}

\subsection{The model space of $Pol_2(n)$}
Denote by $Pol_2(n;\lambda)$ the space of $n$-gons (possibly, with self-intersections) of the total length $2\lambda$ on the plane.  In the case $\lambda=1$ we will write $Pol_2(n)$ for short. This space consists of collections of $n$ vectors 
$(\vec{e}_1, \dots, \vec{e}_n) \in (\R^2)^n$, which satisfy Conditions --1 and --2 (see Introduction). 
Denote by $V_2(\R^n)$ the Stiefel manifold of all pairs of orthonormal vectors in $\R^n$. 
A point  $(Z_1,Z_2) \in V_2(\R^n)$ admits coordinates $Z_j=(z_{1,j}, \dots, z_{n,j})$, $z_{i,j} \in \R$, $j=1,2$, which satisfy the properties:

--a.  $\sum_{i=1}^n (z_{i,1},z_{i,2}) = 0$, $1 \le i \le n$.

--b.  $\sum_{i=1}^n z_{i,j}^2 = 1$, $j=1,2$.

Define the quadratic mapping
$$H: V_2(\R^n) \to Pol_2(n) $$
\begin{eqnarray}\label{H}
H(Z_1,Z_2) = ((z_{1,1}^2-z_{1,2}^2,2z_{1,1}z_{1,2}), \dots ,(z_{n,1}^2-z_{n,2}^2,2z_{n,1}z_{n,2})).
\end{eqnarray}
One may extend this mapping to the space $Pol_2(n;\lambda)$, obviously.

It is convenient to write the coordinates $(z_{i,1},z_{i,2})$ as a point at the imaginary complex  plane in the quaternion space as follows:
$\z_i = z_{i,1} + \j z_{i,2} \in \C \subset \H = \C \oplus \j \C$, $\j^2=-1$. Then for the mapping $H$ we get: 
$$H(\Z)=(\z_1^2,\dots,\z_n^2),$$
$\Z = (\z_1, \dots, \z_n)$. 

Let us check that  the image of $H$ belongs to the subspace $Pol_2(n) \subset (\R^2)^n$.
The condition --1 (see Introduction) is deduced from the conditions --a and --b:
$$ \sum_{i=1}^n \vec{e}_i = \sum_{i=1}^n \z_i^2  = (\sum_{i=1}^n (z_{i,1}^2-z_{i,2}^2), 2\sum_{i=1}^n z_{i,1}z_{i,2})=(0,0).$$

The condition  --2 is proved  as following:
$$ \sum_{i=1}^n \vert \vec{e}_i \vert = \sum_{i=1}^n (\z_i \bar{\z}_i)^2 = \sqrt{\sum_{i=1}^n (z^2_{i,1} - z^2_{i,2})^2 + 4z_{i,1}z_{i,2}}.$$
using --a and --b, we get:
$$\sum_{i=1}^n \vert \vec{e}_i \vert  = \sqrt{\sum_{i=1}^n (z_{i,1}^2 + z_{i,2}^2)^2} = 2.$$

It is easy to see that $H$ is not a homeomorphism, this is a ramified covering of the multiplicity $2^n$.
The space $V_2(\R^n)$ admits the standard fibration:
$$ S^{n-2}  \subset V_2(\R^n) \mapsto S^{n-1}, $$
the standard metric in the image and in the fiber of this fibration induces a metric on $V_2(\R^n)$. This metric determines the probability measure on the space $Pol_2(n)$: the measure of a subset $V \subset Pol_2(n)$ is well-defined as the 
volume of its preimage $H^{-1}(V) \subset V_2(\R^n)$.

\subsection{The model space of $Pol_3(n)$}

The model space of $Pol_3(n)$ is defined by a "complexification"' of the model space $Pol_2(n)$.
Let us consider the following mapping, which is called the Hopf mapping:
\begin{eqnarray}\label{Hopf}
Hopf: \H  \mapsto \R^3, \quad Hopf(\q) = \bar{\q}\i \q,
\end{eqnarray}
where  $\q = q_0 + q_1\i +q_2\j + q_3 \k \in \H$ is a quaternion, $\bar{\q} = q_0 - q_1\i - q_2\j - q_3 \k$, where $\i,\j,k, \in \H$ are the standard generators.

The restriction of the mapping $Hopf$ to the plane $(q_0 + q_2 \j)$ is given by the formula:
$$Hopf(q_0 + q_2 \j) = (q_0^2 - q_2^2, 0, 2q_0q_2).$$
This mapping coincides with the mapping 
$$q_0 + q_2 \j \mapsto (q_0 + q_2 \j)^2,$$ 
where the image is identified with the plane $(z_1,0,z_2 \j) \subset \R^3$.

The restriction of the mapping $Hopf$ to the plane $(q_0 + q_3 \k)$ is given by the formula:
$$Hopf(q_0 + q_3 \k) = (q_0^2 - q_3^2, -2q_0q_3, 0).$$
This mapping coincides with the mapping 
$$q_0 + q_3 \k \mapsto (q_0 + q_2 \k)^2,$$ where the image is identified with the plane $(z_1,-z_2 \k,0) \subset \R^3$.

Denote by $V_2(\C^n)$ the complex Stiefel manifold of all pairs of Hermitian orthonormal vectors in $\C^n$. 
A point  $(Q_{0,1},Q_{2,3}) \in V_2(\C^n)$ admits coordinates $Q_{0,1}=(q_{1,0} + \i q_{1,1}, \dots, q_{n,0}+ \i q_{n,1})$, 
$Q_{2,3}=(q_{1,2} + \i q_{1,3}, \dots, q_{n,2}+ \i q_{n,3})$, 
where $q_{i,1}=q_{i,0}+\i q_{i,1} \in \C$, which satisfy the properties: 

--a'.  $\sum_{i=1}^n (q_{i,0}+\i q_{i,1})(q_{i,2} - \i q_{i,3}) = 0$, $1 \le i \le n$.

--b'.  $\sum_{i=1}^n q^2_{i,0} + q^2_{i,1} = 1$, $\sum_{i=1}^n q_{i,2}^2 + q_{i,3}^2 = 1$.

In \cite{C-G-K-Sh} by a straightforward calculation is proved that 
the coordinatewise extension of the mapping $(\ref{Hopf})$ defines the required mapping:
\begin{eqnarray}\label{Hopff}
Hopf: V_2(\C^n) \to Pol_3(n)
\end{eqnarray}
$$H(Q_0,Q_1,Q_2,Q_3) = (q_{1,0}^2+q_{1,1}^2-q_{1,2}^2-q_{1,3}^2, \quad 2q_{1,1}q_{1,2}- 2q_{1,0}q_{1,3},\quad 2q_{1,0}q_{1,2} + 2q_{1,1}q_{1,3}, \dots ,$$
$$q_{n,0}^2+q_{n,1}^2-q_{n,2}^2-q_{n,3}^2,\quad 2q_{n,1}q_{n,2}- 2q_{n,0}q_{n,3},\quad  2q_{n,0}q_{n,2} + 2q_{n,1}q_{n,3}).$$

\section{Proof of the main result}

\subsection{The normalization of preimages of an $n$-gone in the model space}

The restriction of the mapping $(\ref{Hopf})$ to the sphere $S^3 \subset \H$ of the radius $\sqrt{2}$
(as well as to a sphere of an arbitrary radius) is invariant with respect to the (-right, or -left) multiplication on a complex number $\cos(\theta) + \i \sin(\theta)$. For the left multiplication we get:
$$\overline{(\cos(\theta) + \i \sin(\theta) \q)}\i  (\cos(\theta) + \i \sin(\theta))\q = $$
$$\bar{\q} (\cos(\theta) - \i \sin(\theta)) \i (\cos(\theta) + \i \sin(\theta)) \q = \bar {q} \i \q.$$
For the right multiplication the proof is analogous.

We shall call a pair $(Q_{0,1},Q_{2,3})$ is normalized, if   $q_{i,1}=0$ for an arbitrary $i$, $1 \le i \le n$.
For a normalized pair one could write $Q_{0,1}=(Q_0,0)$.
Let us prove that an arbitrary pair $(Q_{0,1},Q_{2,3})$, which satisfy properties --a', --b', admits the normalization.

Define a complex number $\alpha_i$, $\vert \alpha_i\vert=1$, by the following formula:
$$\alpha_i = \frac{\sqrt{q_{i,0}^2 + q_{i,1}^2}}{q_{i,0} + \i q_{i,1}},$$
if $q_{i,1} \ne 0$, and $\alpha_i = \sign(q_{i,2})$, if $q_{i,1}=0$.
Define the normalization 
$(Q_{0,1},Q_{2,3}) \mapsto (Q'_{0,1}=(Q'_0,0),Q'_{2,3})$ by the following formula:
$$Q'_{0,1}=(q_{1,0} + \i q_{1,1})\alpha_1, \dots (q_{n,0} +\i  q_{n,1})\alpha_n,$$
$$Q'_{2,3}=(q_{1,2} + \i q_{1,3})\alpha_1, \dots (q_{n,2} +\i  q_{n,3})\alpha_n.$$

The images by the  mapping $(\ref{Hopff})$ of the pair of vectors
$(Q_{0,1},Q_{2,3})$ and of its  normalization $(Q'_0,Q'_{2,3})$ coincide. For the normalized pair $(Q'_0,Q'_{2,3})$
the formula of the mapping $(\ref{Hopff})$ is much simple. We get:
\begin{eqnarray}\label{Hopfnorm}
\begin{array}c
H(Q'_0,Q'_{2,3}) = (q_{1,0}^2-q_{1,2}^2-q_{1,3}^2, \quad  - 2q_{1,0}q_{1,3},\quad 2q_{1,0}q_{1,2}, \dots , \\
q_{n,0}^2-q_{n,2}^2-q_{n,3}^2,\quad - 2q_{n,0}q_{n,3}, \quad  2q_{n,0}q_{n,2}).
\end{array}
\end{eqnarray}

Let us prove that the formula $(\ref{Hopfnorm})$ is defined as a superposition of two formulas $(\ref{H})$.
More precisely, a point in $Pol_3(n)$, which parametrizes  an $n$-gon $L \subset \R^3$ of the length $2$, is completely reconstructed from the two projections $L_{x,y} \subset \R^2(x,y,0) \subset \R^3$, 
$L_{x,z} \subset \R^2(x,0,z) \subset \R^3$ onto the standard planes as follows. 

Consider the projection $L_{x,z}$ of $L$ and denote by $l_{x,z}$ the
total length of this projection, $l_{x,y} \le 2$. The projection $L_{x,z}$ is a collection of vectors 
$\vec{e}_i$, which satisfy Condition --a, and a modification of Condition --b, where in the right side of the first and the second formulas we get the constant $l_{x,y}$ instead of $1$. Take a preimage $(Q'_0,Q'_2)$ of $L_{x,z}$ with respect to the mapping
$H$, $H(Q'_0,Q'_2)=L_{x,z}$, $Q'_0=(q'_{1,0}, \dots, q'_{n,0})$, $Q'_2=(q'_{1,2}, \dots q'_{n,2})$.  Let us fix 
a pair $(Q'_0,Q'_2)$ by the condition $Q'_{i,0} \ge 0$, $1 \le i \le n$ (in the case $Q'_{i,0} =0$ let us agree that $Q'_{i,2} \ge 0$ to have a well-defined  lift of $H$). Each vector of the pair $(Q'_1,Q'_2)$ belongs to the 
sphere of the radius $\sqrt{l_{x,y}}$ in $\R^n$.

Apply the same construction to the projection $L_{x,z}$ of $L$. Denote by $l_{x,z}$ the
total length of this projection, $l_{x,z} \le 2$. The following equation, obviously, is satisfied:
$$ l_{x,y} + l_{x,z} - \lambda_0 = 2, $$
where 
$$\lambda_0 = \sum_{i=1}^n (q')^2_{i,0}. $$
Define the vector $Q'_3$, $Q'_3=(q'_{1,3}, \dots q'_{n,3})$, such that for the pair of vector $(Q'_0,Q'_3)$
the equation $H^{conj}(Q'_0,Q'_3)=L_{x,y}$, where   $H^{conj}$ is the composition of $H$ with the reflexion $(x,y) \mapsto (x,-y)$, is satisfied.

The triple $(Q'_0,0,Q'_2,Q'_3)$ of vectors determines the normalized preimage of  $L$.
For an arbitrary preimage $(Q_0,Q_1,Q_2,Q_4)=Hopf^{-1}(L)$ the normalized preimage  $(Q'_0,0,Q'_2,Q'_3)$ is well-defined.  However, the correspondence $L \mapsto  (Q'_0,0,Q'_2,Q'_3)$ is not continuous.

\subsection{A construction with  Cayley numbers related with $Pol_2(7)$}

In this subsection we shall prove the following result, which is, probably, interesting by itself.

\begin{lemma}\label{lemma}
There exists a one-parameter isometric isotopy of the space $Pol_2(7;\lambda)$ to itself, which joins each closed polygon of the length $\lambda$ in the standard plane with its mirror image with the axis $\{x=0\}$. Moreover, this isotopy keeps the projection of each polygon onto the axis $\{x=0\}$.
\end{lemma}

\subsubsection*{Proof of Lemma $\ref{lemma}$}
Assume, for simplicity, that $\lambda =1$, the general case is analogous.
Let us assume that the each vector in the pair $(Z_1,Z_2)$, which determines a point of $V_2(\R^7)$ is represented as 
an imaginary Cally unit, $\vert\vert Z_1 \vert \vert = \vert\vert Z_2 \vert \vert = 1$:
$$Z_j = z_{j,1} \i + z_{j,2} \j + z_{j,3} \k + z_{j,4} \e + z_{j,5} \f + z_{j,6} \g + z_{j,7} \h, \quad j=1,2. $$
Then the vector $Z_1$ determines a one-parameter subgroups in the Cally units $S^7 \subset \Ca$, which contains the elements $\{1,Z_1\}$. Let $Z(\theta)$, $0 \le \theta \le 2\pi$, is the natural parameter of elements of this subgroup, such that 
$Z(\frac{\pi}{2})=Z_1$. Define the following $\theta$--family $Y(\theta)$ of Cally numbers by the formula:
$$ Y(\theta) = Z_{2} \cdot Z(\theta).$$
Obviously, we have $Y(0)=Z_2$,  $Y(\pi) = -Z_2$. Moreover, $Y(\theta)$ is an imaginary Cally number for an arbitrary $\theta$. To proof this, recall that the vector $Z_2$  is orthogonal to $Z^{-1}(\theta)$ for an arbitrary $\theta$.   
The Cally product is represented by an orthogonal transformation of $\R^8$, therefore, 
$Z_2 \cdot Z(\theta)$ is orthogonal to $Z^{-1}(\theta) \cdot Z(\theta) = 1$. By the same argument $Y(\theta)$ is orthogonal to $Z_1$.

The pair $(Z_1,Y(\theta))$ represents a path in $V_2(\R^7)$, which joins the point $(Z_1,Z_2)$ with the point
$(Z_1,-Z_2)$. The projection of the polygon $H(Z_1,Y(\theta))$ onto the line $\{x=0\}$ keeps fixed, because this projection is completely determined by the vector $Z_1 \in \R^7$, and this vector is not changed along the path. But $H(Z_1,-Z_2)$ represents a 7-gone, which is obtained from $H(Z_1,Z_2)$ by the reflection with respect to the line $\{ x=0\} \subset \R^2(x,y)$.

Let us complete the proof. Take an arbitrary 7-gone $L \subset \R^2$, $L \in Pol_2(7)$, which is represented by a collection of coordinates,  satisfied Conditions --a and --b. Take an arbitrary preimage $(Z_1,Z_2) \in H^{-1}(L) \in V_2(\R^7)$. A path $Z_1,Y(\theta))$ is well-defined. To prove that the path $(Z_1,Y(\theta))$ determines an isometric
isotopy of $Pol_2(7)$ to itself, it is sufficient to proof that the projection $H(Z_1,Y(\theta))$ of the path depends no of a lift  $(Z_1,Z_2)$ of $L$. 

The subset $H^{-1}(L) \in V_2(\R^7)$ contains $2^7$ pairs of vectors (probably, less, if a $7$-gone $L$ is degenerated).  An arbitrary two vectors of this collection is related by an automorphism $\alpha \in Aut(\Ca)$, which sends
a subset of the standard generators of $\Ca$ to the antipodal, and all the last standard generators sends to itself.
Obviously, because $H$ is quadratic over each coordinate, by the formula $(\ref{H})$ for an arbitrary $\alpha$ we get: $\alpha \circ H = H$. Lemma $\ref{lemma}$ is proved.

\subsubsection*{Proof of Theorem $\ref{main}$}

Let $L \in Pol_3(7)$ a 7-gons,  $(Q_0,Q_1,Q_2,Q_3) \in V_2(\C^7)$ be a preimage of $L$, $Hopf((Q_0,Q_1,Q_2,Q_3))=L$.
Consider the normalization $(Q'_0,0,Q'_2,Q'_3)$ of this preimage. Apply Lemma \ref{lemma} twice: for the projection $L_{x,z}$ with the preimage
$(Q'_0,Q'_2)$, $H(Q'_0,Q'_2)=L_{x,z}$, $L_{x,z} \in Pol_2(7;l_{x,z})$, and for the projection $L_{x,y}$ with the preimages
$(Q'_0,Q'_3)$, $H^{conj}(Q'_0,Q'_3)=L_{x,y}$, $L_{x,y} \in Pol_2(7;l_{x,y})$.  
 The homotopies for  the polygons $(Q'_0,Q'_2)$, $(Q'_0,Q'_3)$ denote by $L_{x,z}(\theta)$,  $L_{x,y}(\theta)$ correspondingly. By construction $L_{x,z}(0)=L_{x,z}$, $L_{x,z}(\pi)=-L_{x,z}$, $L_{x,y}(0)=L_{x,y}$, $L_{x,z}(\pi)=-L_{x,y}$, where $-L_{x,z}$, $-L_{x,y}$ are the mirror images of $L_{x,z}$, $L_{x,y}$ with respect to the line $\{x=0\}$, correspondingly. Moreover, the projections of $L_{x,z}(\theta)$ and of $L_{x,y}(\theta)$ onto the axis $\{x=0\}$ coincide and depend no of $\theta$. Therefore for an arbitrary $L \in Pol_3(7)$ the required homotopy $L(\theta)$ from $L=L(0)$ to its mirror image $-L=L(\pi)$ is well-defined. 

To prove that the family  of homotopies $Y(\theta)$ determines an isotopy $\Theta$ of the space $Pol_3(7)$ into itself, it is sufficiently to prove that the homotopies continuously depend on $L$. This follows from the fact, that $L(\theta)$ depends no of a lift
$Hopf^{-1}(L)$. By Lemma $\ref{lemma}$ the homotopy of the each projections depends no of its lifts $H^{-1}$. 
The construction of the homotopies $L(\theta)_{x,z}$, $L(\theta)_{x,y}$ depends continuously of its initial points
$H^{-1}(L(\theta)_{x,z})$, $(H^{conj})^{-1}(L(\theta)_{x,y})$. 
Then, even if for two closed polygons $L_1, L_2 \in Pol_3(7)$ the reduced vectors of the preimages jump, the projections of reduced vectors coincide with $L_1$, $L_2$ correspondingly, and $L(\theta)$ continuously depends of $\theta$.

The homotopy $\Theta$ is an isometric isotopy  by the construction, evidently, this homotopy is a bijection. Theorem $\ref{main}$ is proved.

\[  \]
\[  \]
Moscow, Troitsk, IZMIRAN $\qquad \qquad \qquad \qquad$

pmakhmet@mi.ras.ru  $\qquad \qquad \qquad \qquad$


\begin{thebibliography}{99911122}

\bibitem [C-G-K-Sh] {C-G-K-Sh}
Jason Cantarella, Alexander Y. Grosberg, Robert Kusner, and Clayton Shonkwiler, {\it The Expected Total Curvature of Random Polygons}, arXiv:1210.6537v1 [math.DG] 24 Oct 2012.


\bibitem[J]{J}
Ioan James, {\it The topology of Stiefel manifolds} Oxford, (1976)
\end{thebibliography}
\end{document}